\documentclass[12pt,reqno]{amsart}

\PassOptionsToPackage{hyphens}{url}\usepackage[
  hypertexnames=false,
  hidelinks=true
]{hyperref}
\usepackage{amssymb}
\usepackage{fullpage}
\usepackage{amsthm}
\usepackage[noabbrev,nameinlink]{cleveref}
\usepackage[shortlabels]{enumitem}
\usepackage[british]{babel}

\usepackage[
  backend=biber,
  style=ext-numeric,
  natbib=true,
  maxbibnames=99,
  giveninits=true,
  url=false,
  isbn=false,
  doi=false,
  articlein=false,
  mincitenames=99,
  maxcitenames=100,
  urldate=long,
  sorting=nyt
]{biblatex}
\DeclareFieldFormat{pages}{#1}

\newtheorem{thm}{Theorem}[section]
\newtheorem{lemma}[thm]{Lemma}

\crefname{thm}{theorem}{theorems}
\Crefname{thm}{Theorem}{Theorems}

\newcommand{\PP}{\mathbb{P}}
\newcommand{\NN}{\mathbb{N}}

\newcommand{\FF}{\mathbb{F}}
\newcommand{\cal}[1]{\mathcal{#1}}
\newcommand{\eps}{\varepsilon}

\title{
  Szemer\'edi's theorem with random differences \\
  and Arithmetic Kakeya
}

\author{Gabriel Dahia}

\address{Institute of Mathematics, EPFL, Lausanne, Switzerland}
\email{gabriel.dahiafernandes@epfl.ch}

\begin{document}

\begin{abstract}
  We prove that there exists an absolute constant $c > 0$ such that the
  following holds for every $k \ge 3$ and $p \ge n^{-c / \log \log k}$.
  Letting $\delta > 0$ and $D \sim [n]_p$, with high probability every
  ${S \subset [n]}$ with $|S| \ge \delta n$ contains a $k$-term arithmetic
  progression whose common difference lies in $D$.
  This improves the best known bounds for the threshold of Szemer\'edi's theorem
  with random differences for all sufficiently large $k$, previously of the form
  ${p \ge n^{-2/k+o(1)}}$, due to \citeauthor{BC2024}.
  We obtain our improvement via a small cover lemma, in the sense of
  \citeauthor{AP2025+}'s work on random Cayley graphs, and its proof crucially
  relies on a construction of \citeauthor{GR2019} in the context of the
  arithmetic Kakeya conjecture.
  The usage of the latter was suggested to us by GPT 6 Astra.
\end{abstract}

\maketitle

\section{Introduction}

In \citeyear{S1975}, \citet{S1975} proved that subsets of $[n]$ that have no
$k$-term arithmetic progressions have size at most $o(n)$.
Since then, there have been many important generalisations of this result.
Particularly relevant here is the work of \citet{BL1996}, who established that
the same holds if we require the subsets of $[n]$ to avoid only $k$-APs whose
common difference is a perfect $\ell$-th power, for some fixed $\ell \in \NN$.
The set of common differences in this case is very sparse, but also very
structured, so a natural question is how sparse a typical set of common
differences can be while retaining that property.

One formulation of this problem, known as \citeauthor{S1975}'s theorem with
random differences, is the following.
Determine the threshold $p = p(n)$ such that, if $D$ is a $p$-random subset of
$[n]$, then with high probability every subset of $[n]$ that avoids $k$-APs with
common difference in $D$ has size $o(n)$.
The first non-trivial bounds for this problem were established independently by
\citet{C2011} and \citet{FLW2012}, who proved that this holds whenever
$p \ge n^{-1 / 2^{k-2} +o(1)}$.
This was later improved for $k \ge 5$ to $p \ge n^{-1/k+o(1)}$ by
\citet{FLW2016}.
The best current bound is that $p \ge n^{-2/k+o(1)}$ suffices, first established
by \citet{BDG2019} when $k$ is even (see also \cite{BG2020}) and then by
\citet{BC2024} when it is odd.

\citet{FLW2016} conjectured, however, that these bounds are still far from the
truth: they believe that $p n / \log n \to \infty$ should suffice for fixed $k$.
\citet{A2020} moreover proved that if their conjecture is true, then it depends
on the structure of the integers: we now know (see also~\cite{B2021, AT2026+})
that the analogous question in the setting of $\FF_q^d$ for a fixed prime $q \ge k$
requires $p \gtrsim (\log n)^{k - 1} / n$, where $n = q^d$.

Our main result is an improvement over the $p \ge n^{-2/k+o(1)}$ bound for all
sufficiently large $k$.
To state it, let $[n]_p$ denote a $p$-random subset of $[n]$ and
\begin{equation*}
  \cal{D}_k(S) = \big\{d \in \NN : \exists \{a, a + d, \ldots, a + (k - 1)d\} \subset S \big\}
\end{equation*}
denote the set of common differences of the $k$-APs contained in $S$.

\begin{thm}
  \label{stmt:main}
  There exists $c > 0$ such that the following holds for every $k \ge 3$ and
  $\delta > 0$.
  Let $C = C(\delta, k) > 0$ and $n \in \NN$ be sufficiently large.
  If $p \ge C n^{-c / \log \log k}$, then
  \[
    \PP\big(\exists \, S \subset [n], \, |S| \ge \delta n : \cal{D}_k(S) \cap [n]_p = \emptyset \big) \to 0
  \]
  as $n \to \infty.$
\end{thm}

While the strategies of both \citeauthor{BDG2019} and \citeauthor{BC2024} are
technical, relying on connections to coding theory, our approach is more
elementary.
The main inspiration for our work comes from the problem of determining the
independence number of sparse random Cayley sum
graphs~\cite{CDM2024, N2025+, AP2025+}, which is closely related to the
(non-trivial) case of \citeauthor{S1975}'s theorem with random differences for
$k = 2$.

A common feature of the strategies in \cite{CDM2024, N2025+, AP2025+}
(inspired by previous work of \citet{A2007, A2013}) is to replace a
union bound over all candidate independent sets by one over a much
smaller family.
In the formulation of \citet{AP2025+}, this family is a small cover
of sumsets: for every large set $A$, its sumset $A+A$ contains a
set in the cover.
The sets in the cover are moreover required to be sufficiently large.
These two properties, combined with the fact that the family is small, imply
that taking a union bound over the cover succeeds, whereas a naive union bound
over all independent sets fails.

To understand the connection between the approach in the Cayley context and the
current problem, let $\cal{H}_k$ be the hypergraph whose vertex set is $[n]$ and
each hyperedge corresponds to a non-trivial $k$-AP.
Further let $\cal{H}_{k, p}$ be the random hypergraph obtained by keeping only
the edges of $\cal{H}_k$ that have common difference in $[n]_p$.
In this setting\footnote{Indeed, generalising the random sampling argument of
\cite[Lemma~2.7]{CFPY2024+} (see also \cite[Theorem 2.1]{CDM2024}) to
$\cal{H}_{k, p}$ results in a small cover of size $\exp(n^{1-1/k+o(1)})$, which
in turn implies a (suboptimal) bound of $p \ge n^{-1/k+o(1)}$.},
\Cref{stmt:main} is bounding the probability that
$\alpha(\cal{H}_{k, p}) \ge \delta n$.

We similarly deduce \Cref{stmt:main} from a small cover
\namecref{stmt:fingerprints}, \Cref{stmt:fingerprints} below.
The main difference in our argument is that instead of each set $F$ in the cover
being sufficiently large, it will have large $\cal{D}_k(F)$.

\begin{lemma}
  \label{stmt:fingerprints}
  There exists $c > 0$ such that the following holds for every $k \ge 3$ and
  $\delta > 0$.
  Let $C = C(k) > 0$ and $\gamma = \gamma(\delta, k) > 0$, and let $n$ be
  sufficiently large.
  There exists a collection $\cal{F} \subset 2^{[n]}$ with the following
  properties:
  \begin{enumerate}[(a)]
    \item $|\cal{F}| \le C n \exp(n^{1-c/\log \log k})$.\label{item:sizeOfFamily}
    \item For each $S \subset [n]$ with $|S| \ge \delta n$, there exists
          $F \in \cal{F}$ such that
          $F \subset S$.\label{item:fingerprintForEverySet}
    \item Every $F \in \cal{F}$ satisfies $|\cal{D}_k(F)| \ge \gamma n$.
          \label{item:diffsInF}
  \end{enumerate}
\end{lemma}

To see why \Cref{stmt:fingerprints} is useful, we give the short deduction of
\Cref{stmt:main}.

\begin{proof}[Proof of \Cref{stmt:main}]
  Let $\cal{F}$ be given by \Cref{stmt:fingerprints} for $k$, $\delta$ and $n$.
  By \Cref{stmt:fingerprints}\ref{item:fingerprintForEverySet}, observing that
  $\cal{D}_k(F) \subset \cal{D}_k(S)$ holds trivially when $F \subset S$, we
  have
  \begin{equation}
    \label{eq:unionBound}
    \PP(\exists \, S \subset [n], \, |S| \ge \delta n : \cal{D}_k(S) \cap [n]_p = \emptyset) \le \sum_{F \in \cal{F}} \PP(\cal{D}_k(F) \cap [n]_p = \emptyset)
  \end{equation}
  by taking a union bound over all $F \in \cal{F}$.
  Moreover, it follows from \Cref{stmt:fingerprints}\ref{item:diffsInF} that
  \begin{equation}
    \label{eq:probabilityThatMisses}
    \PP\big(\cal{D}_k(F) \cap [n]_p = \emptyset\big) \le (1 - p)^{|\cal{D}_k(F)|}
                                                     \le \exp(- \gamma p n),
  \end{equation}
  for every $F \in \cal{F}$.
  Substituting \eqref{eq:probabilityThatMisses} and
  \Cref{stmt:fingerprints}\ref{item:sizeOfFamily} in \eqref{eq:unionBound}, we
  conclude that
  \begin{equation*}
    \sum_{F \in \cal{F}} \PP\big(\cal{D}_k(F) \cap [n]_p = \emptyset\big)
    \le \exp\!\left( n^{1-c / \log \log k} + 2 \log n - \gamma p n \right)
    \le \exp\!\left( - \gamma n^{1-c / \log \log k} \right)
  \end{equation*}
  by our choice of $p \ge C n^{-c / \log \log k}$, for any
  ${C > 3 \gamma^{-1}}$.
\end{proof}

\section{Proof of \Cref{stmt:fingerprints}}

In the proof of \Cref{stmt:fingerprints}, we will use two results.
The first is a standard result of \citet{V1959}, a ``supersaturation'' form of
\citeauthor{S1975}'s theorem.

\begin{thm}[\cite{V1959}, see {\cite[Theorem 9.3.1]{Z2023}}]
  \label{stmt:varnavides}
  Let $k \in \NN$, $\delta > 0$ and let $n \in \NN$ be sufficiently large.
  If $S \subset [n]$ satisfies $|S| \ge \delta n$, then it contains at least
  $\eps n^2$ $k$-APs, for $\eps = \eps(\delta, k) > 0$.
\end{thm}

The other, and crucial, input is the following result of \citet{GR2019}, which
they relate to the arithmetic Kakeya conjecture of Katz and Tao.

\begin{thm}[{\cite[Theorem 1.7]{GR2019}}]
  \label{stmt:greenRuzsaArithmeticKakeya}
  There exists $c > 0$ such that the following holds for every $k \ge 3$ for
  some $C = C(k) > 0$.
  For every sufficiently large $n \in \NN$, there exists $Y \subset [C n]$ with
  $|Y| \le n^{1-c / \log \log k}$ which contains, for every $d \in [n]$, a
  $k$-AP with common difference $d$.
\end{thm}

The construction of $Y$ takes a single page of \citeauthor{GR2019}'s
paper~\cite[Section 5]{GR2019}.
They first construct a smaller set $Y' \subset [C_k]$, selecting the starting
point of each $k$-AP to be the square of its common difference modulo $Q$,
the product of the first few odd primes.
The Chinese remainder theorem allows them to bound the size of the set of the
$j$-th terms of all the $k$-APs for each $j \in [k]$, and then of the entire
$Y'$ by $Q^{1-c / \log \log k}$.
They finally obtain $Y$ from $Y'$ using the base $Q$ representations of the
common differences $d \in [n]$.

As will become clear in the proof of \Cref{stmt:fingerprints} below, improving
the upper bound on $|Y|$ in \Cref{stmt:greenRuzsaArithmeticKakeya} would give a
corresponding improvement for the threshold in \citeauthor{S1975}'s theorem with
random differences.
This reduction then connects two problems in \citeauthor{G2025++}'s list of 100
open problems, which also highlighted the case $k = 3$ of \citeauthor{S1975}'s
theorem with random differences~\cite[Problems 10 and 37]{G2025++}.

However, it is worth noting that no improvement in
\Cref{stmt:greenRuzsaArithmeticKakeya} would lead to a polynomial improvement
for the known bounds when $k \in \{3, 4\}$.
To see why, observe that the case $k = 2$ is equivalent to $[n] \subset Y - Y$,
and that forces $|Y| \ge n^{1/2}$.
On the other hand, the results of \citeauthor{BC2024} (in fact, even those in
\cite{BDG2019, BG2020} for $k = 4$) already imply that the threshold in these
cases is at most $n^{-1/2+o(1)}$.

\begin{proof}[Proof of \Cref{stmt:fingerprints}]
  Let $Y$ be the set given by \Cref{stmt:greenRuzsaArithmeticKakeya}, and, for
  each $d \in [n]$, let $P_d \subset Y$ be the $k$-AP with common difference
  $d$, and fix a set $S \subset [n]$ with $|S| \ge \delta n$.
  Observing that each ${P \subset S}$, where $P$ is a $k$-AP, can be written as
  $P = t + P_d$ for some pair ${(t, d) \in \{-C n, \ldots, C n\} \times [n]}$,
  it follows from \Cref{stmt:varnavides} that there are $\eps n^2$ pairs
  $(t, d)$ where $t + P_d \subset S$.
  This implies that there is one $t_*$ such that, letting
  \begin{equation}
    \label{eq:formOfTheFingerprint}
    F = (t_* + Y) \cap S,
  \end{equation}
  we have that the number of $t_* + P_d \subset F$ is at least
  $$
    \frac{\eps n^2}{2 C n + 1} \ge \frac{\eps n}{2C + 1}.
  $$
  Consequently, as every $P_d \subset Y$ has a distinct common difference $d$,
  we conclude that
  \begin{equation}
    \label{eq:propertyOfFingerprint}
    |\cal{D}_k(F)| \ge \frac{\eps n}{2C + 1} = \gamma n
  \end{equation}
  for $\gamma = \eps / (2 C + 1)$.
  To complete the proof, we let $\cal{F}$ be the family of the $F$ satisfying
  \eqref{eq:formOfTheFingerprint} and \eqref{eq:propertyOfFingerprint} for some
  $S \subset [n]$ with $|S| \ge \delta n$, and note that they are fully
  specified by choosing $t_* \in \{-Cn, \ldots, Cn\}$ and one of the
  $2^{|Y|} \le \exp(n^{1-c / \log \log k})$ sets $Y \cap (S - t_*) \subset Y$.
\end{proof}

\subsection*{AI Use}

We prompted GPT 6 Astra Pro to see if the ideas in \cite{CL2026+} could be used
to construct a small cover for the hypergraph encoding $k$-APs.
It replied that it could not use those ideas, but suggested that the
construction from \cite{GR2019} could be used to prove such a result.

\section*{Acknowledgments}

We are grateful to Jo\~ao Pedro Marciano, Rob Morris, Marcelo Campos and Oliver
Janzer for suggestions that improved the paper's presentation.
We also thank them for previous discussions about the problem of
\citeauthor{S1975}'s theorem with random differences.

\def\bibfont{\footnotesize}
\printbibliography

\end{document}